\documentclass[a4paper,notitlepage,twoside,reqno,11pt]{amsart}

\usepackage{bbm,pifont,latexsym,amssymb}
\usepackage{anysize}
\marginsize{3.4cm}{3.4cm}{3cm}{3cm}

\usepackage{dcolumn,indentfirst}
\usepackage[hypertexnames=false,bookmarksopen=true,linktocpage=true,pdfstartview={XYZ null null 1.25}]{hyperref}
\usepackage{amsmath,amssymb,amscd,amsthm,amsfonts,mathrsfs}
\usepackage{color,graphicx,xcolor,graphics,subfigure,extarrows,caption2}
\usepackage{titletoc}

\usepackage{autonum}

\newtheorem{thm}{Theorem}[section]
\newtheorem{lem}[thm]{Lemma}

\newtheorem*{main}{Main Theorem}

\theoremstyle{definition}

\newcommand{\T}{\mathbb{T}}

\newcommand{\Z}{\mathbb{Z}}

\newcommand{\C}{\mathbb{C}}
\newcommand{\N}{\mathbb{N}}
\newcommand{\EC}{\widehat{\mathbb{C}}}

\newcommand{\CV}{\textup{CV}}

\newcommand{\Ext}{\mathcal E}
\newcommand{\Crit}{\operatorname{Crit}}
\newcommand{\Fatou}{\mathcal F}
\newcommand{\Julia}{\mathcal J}

\makeatletter\@addtoreset{equation}{section}\makeatother

\begin{document}

\author[Yuming Fu]{Yuming Fu}
\address{School of Mathematical Sciences, Shenzhen University, Shenzhen 518060, P. R. China}
\email{yumingfuxy@szu.edu.cn}

\author[Weiwei Su]{Weiwei Su}
\address{School of Mathematics, Nanjing University, Nanjing 210093, P. R. China}
\email{suweiwei2024@163.com}

\author[Fei Yang]{Fei Yang}
\address{School of Mathematics, Nanjing University, Nanjing 210093, P. R. China}
\email{yangfei@nju.edu.cn}

%---------------------------------------------------------------------------------------------------------------
\title[Rational maps without Herman rings]{Rational maps with at most three critical values have no Herman rings}

\begin{abstract}
We prove that any rational map with at most three critical values has no Herman rings. This generalizes previous results of Shishikura and Milnor on bi-critical rational maps.
\end{abstract}

% AMS subject classifications (used in AMS journals)

% AMS keywords (used in AMS journals)

% today's date, or fill in whatever date you prefer
\date{\today}

% acknowledge support, etc
% \thanks{This research was partially supported by NSF grant DOA-123456789.}

% dedication
% \dedicatory{Dedicated to Professor}

\maketitle

%----------------------------------------------------------------------------------------------------------------
%{\setcounter{tocdepth}{1}
%\tableofcontents
%}

%----------------------------------------------------------------------------------------------------------------
\section{Introduction}\label{introduction}

Let $f:\EC\to\EC$ be a rational map of degree at least two. Its Fatou set $\Fatou(f)$ is the largest open set on which the iterates $\{f^n\}_{n\in\N}$ of $f$ form a normal
family, and its Julia set is $\Julia(f)=\EC\setminus\Fatou(f)$. A periodic Fatou component $H$ of $f$ is called a \textit{Herman ring} if $H$ is conformally isomorphic to some annulus and if the restriction of $f$, or some iterate of $f$ on $H$, is holomorphically conjugate to an irrational rotation of this annulus.
Unlike other types of periodic Fatou components of rational maps, a Herman ring has no associated periodic point. In general it is difficult to judge whether a given rational map (or a family of rational maps) has a Herman ring or not. There are two known methods for constructing Herman rings. The original method is based on a careful analysis on the analytic diffeomorphism of the circle, due to Herman \cite{Her79}. An alternative method is based on quasiconformal surgery, by pasting together two cycles of Siegel disks, due to Shishikura \cite{Shi87}. Herman's original construction yields period-one Herman rings that are symmetric about a circle, whereas Shishikura's construction is more flexible. In particular, it allows one to construct asymmetric Herman rings and Herman rings of any period.

In \cite{Shi87}, Shishikura established an inequality that bounds the number of Herman rings in terms of the degrees of rational maps, and proved that no quadratic rational map has a Herman ring. In \cite{Mil00b}, Milnor generalized Shishikura's result and proved that no rational map with only two critical points (without counting multiplicity) has a Herman ring.
By the Riemann-Hurwitz formula, a rational map with exactly two critical points is equivalent to that it has exactly two critical values.
In this paper, we generalize Shishikura and Milnor's results to the following.

\begin{main}
Rational maps with at most three critical values have no Herman rings.
\end{main}

This result is sharp in the sense that Herman rings can indeed appear in the cubic rational maps which have \textit{four} critical values (see \cite{Her79}, \cite{Shi87} and \cite[p.\,164]{Mil06}).
Recently some sufficient conditions for rational maps to not admit Herman rings have been established.
It was proved in \cite{Yan17} that if a rational map has at most one infinite critical grand orbit in the Julia set, then this map has no Herman rings.
Hu and Xiao proved a partial result of the Main Theorem: any rational map with exactly three critical values does not have Herman rings of periods $1$, $2$ and $3$.
Their proof is based on the study of Shishikura's puzzle, which in particular implies that the \textit{regularly ramified} rational maps (the local degree is the same at all preimages of any given point) have no Herman rings \cite{HX19}.

In the proof of Hu and Xiao's result \cite{HX19}, they ruled out Herman rings of periods $1$, $2$ and $3$ by comparing the conformal moduli of some annuli. However, to exclude Herman rings of more higher period, one needs to use a conformal invariant to distinguish multiply connected domains in the sense of conformal equivalence.
Lemma \ref{lem:geodesic} in this paper precisely serves this purpose.
The proof of the Main Theorem provides a unified treatment for all periods of Herman rings, as well as for the rational maps with exactly $2$ or $3$ critical values.

One of the motivations for studying the existence of Herman rings is to understand the connectivity of Julia sets. In particular, to prove that the Julia set of a rational map is connected, most of the time the first step is to rule out the existence of Herman rings.
Besides the bi-critical rational maps studied in \cite{Shi87}, \cite{BB99} and \cite{Mil00b}, the Main Theorem can be used to exclude the Herman rings in some special families of rational maps studied in \cite{XQ10}, \cite{DR13}, \cite{XQY14}, \cite{FHK19}.

For some related study of Herman rings of rational maps, see \cite{Shi89}, \cite{BFGH05}, \cite{FH17}, \cite{Yan22} and the references therein.

\medskip
\noindent\textbf{Acknowledgments.}
Part of the proof idea was suggested by ChatGPT.
This work was supported by NSFC (Grant Nos.\,12526526 and 12571093).

\section{Proof of the Main Theorem}\label{sec-proof}

Let $f$ be a rational map of degree at least two. We use $\Crit(f)$ and $\CV(f)$ to denote the sets of critical points and critical values of $f$ respectively (without counting multiplicity).
Suppose $\{H_0,\ldots,H_{p-1}\}$ is a cycle of Herman rings of $f$ of exact period $p\geq 1$, with indices taken modulo $p$.
Then there exist an irrational number $\theta$ and a number $0<r<1$ such that for each $i\in\Z_p$, there is a conformal isomorphism $\phi_i:H_i\to \mathbb{A}_r:=\{\zeta\in\C:r<|\zeta|<1\}$ satisfying
\[
 \phi_i\circ f^p\circ\phi_i^{-1}(\zeta)=e^{2\pi i\theta}\zeta, \text{\quad where } \zeta\in\mathbb{A}_r.
\]
For a given $s\in(r,1)$, we denote $\gamma_0:=\phi_0^{-1}(\T_s)$, where $\T_s:=\{\zeta\in\C:|\zeta|=s\}$ is a round circle with radius $s$.
Set $\gamma_i:=f^i(\gamma_0)$. Then $\{\gamma_i:i\in\Z_p\}$ are pairwise disjoint analytic Jordan curves, essential in their respective Herman rings, and
\begin{equation}\label{eq:cores}
 f(\gamma_i)=\gamma_{i+1} \text{\quad and\quad} f^p(\gamma_0)=\gamma_0.
\end{equation}
We call them a \emph{cycle of core curves}. By varying the initial radius $s\in(r,1)$, we may arrange that all core curves avoid any prescribed finite set.
For a bounded open Jordan disk $D\subset\C$, we denote the exterior of the closure of $D$ by
\[
 \Ext(D):=\EC\setminus\overline D.
\]

\begin{lem}\label{lem:disks}
Suppose $\#\CV(f)\leq 3$ and $\infty\in\Crit(f)$. Let $\{\gamma_0,\ldots,\gamma_{p-1}\}$ be a cycle of core curves avoiding
$\CV(f)$, and let $D_i$ be the bounded Jordan disk enclosed by $\gamma_i$. Then for each $i\in\Z_p$, one of the following two cases happens:
\begin{itemize}
\item $\#(D_{i+1}\cap\CV(f))\leq1$ and $f:D_i\to D_{i+1}$ is conformal; or
\item $\#(D_{i+1}\cap\CV(f))\geq2$ and $f:D_i\to \Ext(D_{i+1})$ is conformal.
\end{itemize}
\end{lem}

\begin{proof}
Let $D_i'$ be the component of $f^{-1}(D_{i+1})$ such that $\gamma_i\subset\partial D_i'$. If $\#(D_{i+1}\cap\CV(f))\leq1$, then $f:D_i'\to D_{i+1}$ is a proper map branched over at most one point. By the Riemann-Hurwitz formula, $D_i'$ is simply connected, and in particular, $D_i'$ is a Jordan disk and $\partial D_i'=\gamma_i$. Since $f:\gamma_i\to\gamma_{i+1}$ has degree one, it follows that  $f:D_i'\to D_{i+1}$ is conformal. Since $\infty\in\Crit(f)$, we conclude that $D_i'=D_i$ and $f:D_i\to D_{i+1}$ is conformal.

If $\#(D_{i+1}\cap\CV(f))\geq2$, then $\#(\Ext(D_{i+1})\cap\CV(f))\leq 1$. By the same argument as above, we have the conformal map $f:D_i\to \Ext(D_{i+1})$.
\end{proof}

\begin{lem}\label{lem:geodesic}
For $n\geq 2$, let
\begin{equation}
\Omega:=\EC\setminus\bigcup_{j=1}^nK_j \text{\quad and\quad} \Omega':=\EC\setminus\bigcup_{j=1}^nK_j'
\end{equation}
be two multiply connected domains, where $K_j$ (resp. $K_j'$) are pairwise disjoint closed Jordan disks.
If $K_{j_0}'\subsetneq K_{j_0}$ for some $j_0$ and $K_j'\subset K_j$ for all $j\neq j_0$, then $\Omega$ and $\Omega'$ are not conformally equivalent.
\end{lem}

\begin{proof}
For each $1\leq j\leq n$, by \cite[Theorem 7.2.5, p.\,129]{KL07}, there exists a unique shortest closed hyperbolic geodesic $\alpha_j$ (resp., $\alpha_j'$) in $\Omega$ (resp., $\Omega'$) which separates $K_j$  (resp., $K_j'$) from all the other complementary components.
Let $\ell_U(\alpha)$ be the hyperbolic length of a rectifiable simple curve $\alpha$ in $U=\Omega$ or $\Omega'$ with respect to the corresponding Poincar\'{e} metric. By the assumption, $\Omega\subsetneq\Omega'$. According to the Schwarz-Pick lemma, this inclusion strictly decreases the Poincar\'e metric. Hence
\[
   \ell_{\Omega'}(\alpha_j')
    < \ell_\Omega(\alpha_j')
    \leq \ell_\Omega(\alpha_j),
   \quad \text{ where } 1\leq j\leq n.
\]
Assume there exists a conformal map $\varphi:\Omega\to\Omega'$. Then $\varphi$ permutes the boundary components and preserve the corresponding hyperbolic lengths of the shortest closed geodesics. Thus, we would have
\[
    \sum_{j=1}^{n}\ell_\Omega(\alpha_j)
    = \sum_{j=1}^{n}\ell_{\Omega'}(\alpha_j')
    < \sum_{j=1}^{n}\ell_\Omega(\alpha_j),
\]
which is a contradiction.
\end{proof}

A subset $A$ of $\EC$ is compactly contained in another subset $B$ if the closure of $A$ is contained in the interior of $B$ and we denote it by $A\Subset B$.

\begin{lem}\label{lem:nested}
Let $\gamma_0,\ldots,\gamma_{p-1}\subset\C$ be a cycle of core curves of a rational map $f$. Let $U$ and $V$ be the bounded disks enclosed by two different core curves, with $U\Subset V$. Then no positive iterate $G=f^k$ restricts to a
conformal isomorphism $G:V\to\Ext(U)$.
\end{lem}

\begin{proof}
Suppose we have a conformal map $G=f^k:V\to\Ext(U)$ for some $k\geq 1$. Note that $G(\partial V)=\partial U$ and $G$ permutes the core curves $\{\gamma_0,\ldots,\gamma_{p-1}\}$.

Start with $B_1=U$. Suppose that $B_1,\ldots,B_n$ have been constructed as bounded Jordan disks with pairwise disjoint closures compactly contained in $V$, and that $G(B_j)=B_{j+1}$ for $1\leq j<n$. Denote
\[
 B_{n+1}:=G(B_n) \text{\quad and \quad}
 \Omega_n:=V\setminus\bigcup_{j=1}^n\overline{B}_j.
\]
Since $G$ is conformal on $V$, $B_{n+1}$ is a spherical Jordan disk bounded by a core curve and $\overline{B}_{n+1}\subset\Ext(U)$. By the injectivity of $G|_V$ and the induction hypothesis, we have
\begin{equation}\label{equ:disjoint}
\overline{B}_{n+1}\cap \overline{B}_j=\emptyset \text{\quad for all } 1\leq j\leq n.
\end{equation}
Then
\begin{equation}\label{eq:punctured-image}
 G(\Omega_n)=\Ext(U)\setminus\bigcup_{j=2}^{n+1}\overline{B}_j.
\end{equation}
There are three possibilities.

\smallskip
\emph{Case 1: $B_{n+1}$ is bounded and $B_{n+1}\Subset V$.}
The induction continues.

\smallskip
\emph{Case 2: $B_{n+1}$ is bounded and is not compactly contained in $V$.}
Its boundary and $\partial V$ are either identical or disjoint. Since $\overline{B}_{n+1}$ is disjoint from $\overline U\subset V$, we have $\overline{B}_{n+1}\Subset\Ext(V)$.
Then \eqref{eq:punctured-image} becomes
\[
 G(\Omega_n)=\EC\setminus\bigcup_{j=1}^{n+1}\overline{B}_j.
\]
The complementary disks of $\Omega_n$ are $\overline{B}_1,\ldots,\overline{B}_n$ and $\EC\setminus V$. The image domain has the same first $n$ complementary disks, while
$\EC\setminus V$ is replaced by the strictly smaller disk $\overline{B}_{n+1}$. This is impossible by Lemma \ref{lem:geodesic}.

\smallskip
\emph{Case 3: $B_{n+1}$ contains infinity.}
Write $B_{n+1}=\Ext(V')$, where $V'$ is the bounded Jordan disk enclosed by its core boundary. By \eqref{equ:disjoint}, we have $\bigcup_{j=1}^n\overline{B}_j\Subset V'$.
Since both $V$ and $V'$ contain $U$ and $\partial V$, $\partial V'$ are core curves, so either $V\Subset V'$, $V'\Subset V$, or $V=V'$. Now by \eqref{eq:punctured-image}, we have
\begin{equation}\label{eq:bounded-image}
 G(\Omega_n)=V'\setminus\bigcup_{j=1}^n\overline{B}_j.
\end{equation}
If $V'\ne V$, the domains $\Omega_n$ and $G(\Omega_n)$ differ by a strict shrinking of exactly one complementary closed disk, which is impossible by Lemma \ref{lem:geodesic}.
If $V'=V$, then $G(\Omega_n)=\Omega_n$ and $\Omega_n$ is contained in the Fatou set of $f$. This is also a contradiction since $\Omega_n$ contains the core curve $\partial U$ and must contain a subset of the Julia set of $f$.

\smallskip
Finally, Case 1 cannot continue indefinitely. Indeed, the period of $\partial V$ under the permutation of core curves induced by $G=f^k$ is $q=p/\gcd(p,k)$. Since $G(\partial V)=\partial U\ne\partial V$, we have $q\geq2$. At every stage, $\partial B_j=G^j(\partial V)$.
If the construction has not stopped earlier, $B_q$ is defined and has boundary $\partial V$.
Its closure is disjoint from $\overline U$, so it cannot be $V$. Hence $B_q=\Ext(V)$, which is the already excluded subcase of Case 3. The proof is complete.
\end{proof}

\begin{proof}[Proof of the Main Theorem]
Let $f$ be a rational map of degree at least two which satisfies $\#\CV(f)\leq 3$.
Assume that $f$ has a $p$-cycle of Herman rings $\{H_i:0\leq i<p\}$ which contains a $p$-cycle of core curves $\{\gamma_i:0\leq i<p\}$.
Conjugate by a M\"obius transformation if necessary, we assume that $\infty\in\Crit(f)$ and $\gamma_i\cap\CV(f)=\emptyset$ for all $i$. Let $D_i$ be the bounded disk bounded by $\gamma_i$ for $0\leq i<p$.

We call a step \emph{interior} if $f(D_i)=D_{i+1}$ and \emph{exterior} if $f(D_i)=\Ext(D_{i+1})$ (these are the only cases by Lemma \ref{lem:disks}). Orient every $\gamma_i$ as the positive boundary of its bounded disk. An interior step preserves these orientations, whereas an exterior step reverses them. The return map on $\gamma_0$ is conjugate to an irrational rotation and preserves orientation. Thus the number $m$ of exterior steps is even.
If $m=0$, then we have $f^p(D_0)=D_0$ and hence $D_0$ is contained in the Fatou set of $f$. This is impossible since $D_0$ contains a subset of the Julia set of $f$.
Thus $m\geq2$.

For each exterior step, its target bounded disk $D_{i+1}$ contains at least two of the three (or two) critical values by Lemma \ref{lem:disks}.
Any two such target disks share a critical value and therefore intersect. Their boundaries are distinct core curves, so these disks form a strictly nested chain.

We list the target disks for all exterior steps in temporal order along the cycle as
\[
 V_1,V_2,\ldots,V_m \text{\quad and\quad} V_{m+1}=V_1.
\]
Starting at $V_j$ and proceeding to the next target disk, all steps are interior until the final exterior step. Hence there exists $k_j\geq1$ such that
\begin{equation}\label{eq:temporal}
 f^{k_j}:V_j\to\Ext(V_{j+1})
\end{equation}
is conformal.
Choose $V_j$ maximal in the nested chain. Its successor $V_{j+1}$ is a different
disk because $m\geq2$, and thus $V_{j+1}\Subset V_j$. Applying Lemma \ref{lem:nested} to $U=V_{j+1}$ and $V=V_j$, we obtain a contradiction.
The proof is complete.
\end{proof}

%%%%%%%%%%%%%%%%%%%%%%%%%%%%%%%%%%%%%%%%%%%%%%%%%%%%%%%%%%%%%%%%%%%%

\bibliographystyle{amsalpha}
\bibliography{E:/Latex-model/Ref1}

\end{document}